\documentclass[conference]{IEEEtran}
\usepackage{cite}

\newcommand\blfootnote[1]{%
\begingroup
\renewcommand\thefootnote{}\footnote{#1}%
\addtocounter{footnote}{-1}%
\endgroup
}
\renewcommand{\footnoterule}{%
\kern -3pt
\hrule width 0.49 \textwidth height 0.5pt
\kern 1pt
}
\ifCLASSINFOpdf
   \usepackage[pdftex]{graphicx}
   \DeclareGraphicsExtensions{.pdf,.jpeg,.png}
\else
   \usepackage[dvips]{graphicx}
   \DeclareGraphicsExtensions{.eps}
\fi
   \usepackage{xcolor}

\usepackage{amsmath}

\usepackage{multirow}
\begin{document}
%
\title{Exact Distribution Optimal Power Flow (D-OPF) Model using Convex Iteration Technique}

\author{\IEEEauthorblockN{Rahul Ranjan Jha and Anamika Dubey}
\IEEEauthorblockA{School of Electrical Engineering and Computer Science \\
Washington State University\\
Pullman, WA}}


%


\maketitle

\begin{abstract}
The distribution optimal power flow (D-OPF) models have gained attention in recent years to optimally operate a centrally-managed distribution grid.  On account of nonconvex formulation that is difficult to solve, several relaxation methods have been proposed; the exactness of the solutions obtained from the relaxed models, however, remain a concern. In this paper, we identify one such problem related to radial distribution feeder where second-order cone program (SOCP) relaxation does not yield a solution that is feasible with respect to the original nonlinear OPF model. Specifically, we formulate an OPF model for PV hosting capacity problem to obtain maximum PV capacity that a feeder can integrate without violating the operating constraints. The SOCP relaxation for this problem yields infeasible solutions. To address this concern, we propose a convex iteration technique to simultaneously achieve optimal and feasible OPF solution for the PV hosting (maximization) problem. The proposed approach minimizes the feasibility gap with respect to original nonlinear constraints at each SOCP iteration. The effectiveness of the proposed approach is validated using IEEE-13 node and IEEE-123 node test systems.  
\end{abstract}

\vspace{-0.4cm}
\blfootnote{\noindent \hspace{-0.265 cm}This work was supported in part by the U.S. Department of Energy under contract DE-AC05-76RL01830.} 

%
\IEEEpeerreviewmaketitle

\vspace{-0.4cm}

\section{Introduction}


With the advancement in smart grid technology and increasing penetration of distributed energy systems (DERs), the electric power distribution system is rapidly transforming to an active and bidirectional network. In a centrally managed distribution grid, an optimal power flow (OPF) solver finds multiple applications related to effective management of all grid resources including but not limited to loss minimization, volt-var optimization, and effective management of DERs \cite{OPFsurvey,OPFsurvey2}. The related literature on OPF models from transmission systems is not directly applicable to distribution grid on account of radial feeders, high R/X ratio, and large-variations in bus voltages. Consequently, several researchers have proposed distribution OPF (D-OPF) formulations. D-OPF models have been largely based on two power flow methods: bus-injection model (BIM) and branch-flow model (BFM). Although bus-injection model applies to general radial/mesh feeders, branch flow method is more suitable for modeling radial distribution feeders. Both BIM and BFM based OPF models are non-convex and difficult to solve. 

To address this concern, several relaxed models have been proposed in literature that deal with the problem of non-convexity either using convex relaxation techniques \cite{SHLowCR1} or using linear approximation methods \cite{approx}. Note that a BIM OPF formulation can be relaxed as semidefinite program (SDP) by dropping the rank-1 constraint \cite{XBai} and a BFM OPF model can be relaxed as a second-order cone program (SOCP) by relaxing the quadratic equality constraints \cite{RAJabr}. Although the proposed SDP and SOCP relaxations result in a convex problem thus reducing the complexity of nonlinear OPF model, the exactness of the solution obtained from relaxed model is still of concern. Consequently, several researchers have attempted to prove the exactness of relaxed OPF models. Sufficient conditions were provided under which the relaxed SOCP and SDP models are exact \cite{DKMolzahn,Low1, Low2,MFarivar,LGan}. Further, it was proved in \cite{Sbose}, that SDP, chordal, and SCOP convex relaxation techniques are equivalent for the radial network while for mesh networks, SDP and chordal relaxation perform better. 

The focus of this work is on D-OPF relaxation techniques applied for radial distribution system. Therefore, we focus our attention to SOCP problems. It has been proved in the existing literature that the SOCP relaxation is exact for the radial distribution feeders under certain conditions \cite{Low2}. Authors also claim that even when these conditions are not exactly satisfied, the relaxed SOCP formulation yields a solution that is optimal and feasible with respect to the original problem. Further, in \cite{Low2}, authors have provided an excellent visualization of conditions for the exactness of relaxed SOCP models. With the help of a two-bus model, it is demonstrated that the SOCP model leads to an optimal solution that lies at the boundary of the second-order power cone, thus achieving feasibility.  

Unfortunately, the exactness of the SOCP relaxation (for radial feeder) is contingent upon the choice of objective function. Although, SOCP relaxation is exact for most minimization problems (in power flow variables), a maximization problem yields infeasible power flow solution. Notice that although most OPF problems do relate to minimizing a cost function of power flow variables, there are relevant cases when a maximization problem needs to be solved. One such case is identifying maximum photovoltaic (PV) penetration limits for the distribution feeder also known as PV hosting problem. This requires maximizing the sum of nodal power injection (from PVs) until system operating constraints (thermal limit, voltage, reverse power flow) are not violated. Solving this optimization problem using SOCP model yield infeasible results that lie inside the second-order cone and not at its boundary. 

To address this gap, the objective of this paper is to propose a convex iteration technique to simultaneously achieve optimal and feasible D-OPF solution for the PV hosting (maximization) problem. Related to this problem, in \cite{iter2} the hosting capacity of the distribution system is obtained, using an iterative SOCP method proposed in \cite{iter1}. The iterative approach detailed in \cite{iter1} is based on generating linear inequality constraints for the relaxed SOCP problem. Finding these linear cuts require solving another optimization problem thus increasing the computational complexity. Alternatively, in this paper, we present an iterative approach that is based on convex iteration sequences. A feasible OPF solution is obtained by solving multiple SOCP iterations of the relaxed OPF problem modeled as SOCP. The iteration sequences are based on linear inequality constraints obtained using a mathematical analysis of quadratic equality constraints as detailed in \cite{Filter}. Note that every iteration solves a SOCP problem where the optimal solutions are updated to make them feasible wrt. the original quadratic equality constraints. The benefit of adding this constraint is that at each iteration, the solution will approach towards the surface of the second-order cone. Thus, the feasibility of the actual problem is guaranteed. We also discuss the reasons for SOCP relaxation not being exact for the maximization problem. Next, with the help of IEEE-13 bus and 123-bus distribution test systems, we analyze the feasibility gap when solving SOCP problem for PV maximization. It is also demonstrated that the exact solutions are obtained by solving multiple iterations of SOCP problem with added linear inequality constraints proposed in this paper.


\section{D-OPF Formulation for PV Hosting Capacity Problem using Branch-Flow Model}
In this section, first, we present the branch-flow model for the distribution system and detail the optimal power flow formulation for the PV hosting capacity problem. Next, we describe the SOCP relaxation by relaxing the quadratic equality constraint and discuss the problem with attaining the exact solution for PV hosting problem (maximization problem) using the relaxed SOCP model. 
\begin{figure}[h]
\vspace{-0.6cm}
\centering
\includegraphics[width=2.5in]{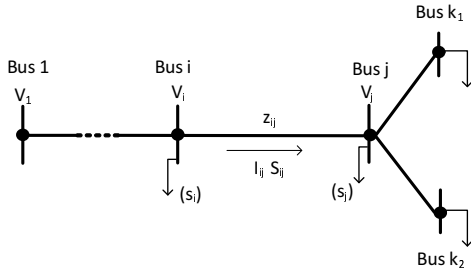}
\vspace{-0.5cm}
	\caption{Topology of radial distribution system }
	\vspace{-0.6cm}
	\label{fig:1}
	\end{figure}
\subsection{Branch Flow Equations}
Consider a directed graph $\mathcal{G}(\mathcal{N, E})$ representing the radial distribution feeder in Fig.1, where, $\mathcal{N}$ is the number of buses and  $\mathcal{E}$ is the number of edges in the graph. The edge ($i,j$) connects nodes $i$ and $j$ where node $i$ is the parent of  node $j$. For each edge $(i,j) \in \mathcal{E}$ associate a complex number $z_{ij} = r_{ij} + \iota x_{ij}$ representing the complex impedance of the line. Also, for each edge ($i,j$), assume that the apparent power flow is $S_{ij} = P_{ij} + \iota Q_{ij}$ and complex line current is $I_{ij}$. For each node $(i) \in \mathcal{N}$, let $V_i$ be the complex voltage and $s_i = p_i + \iota q_i$ be the net apparent power injection (generation minus demand) at corresponding bus $i$.

Then the branch-flow equations for radial distribution feeder, represented by $\mathcal{G}(\mathcal{N,E})$ is given in (1)-(4). Please refer to \cite{MFarivarLow} for further details. 

\vspace{-0.5cm}
\begin{small}
\begin{eqnarray}
P_{ij} &=& \sum_{k:j\to k}P_{jk} + r_{ij}l_{ij} + p_j \hspace{2.2 cm} \forall i \in \mathcal{N}\\
Q_{ij} &=& \sum_{k:j\to k}Q_{jk} + x_{ij}l_{ij} + q_j \hspace{2.2 cm} \forall i \in \mathcal{N}\\
v_j &=& v_i - 2(r_{ij}P_{ij} + x_{ij}Q_{ij}) + (r_{ij}^2 + x_{ij}^2)l_{ij}  \forall i \in \mathcal{N}\\
v_il_{ij} &=& P_{ij}^2 + Q_{ij}^2  \hspace{4 cm} \forall {ij} \in \mathcal{E}
\end{eqnarray}
\end{small}
\vspace{-0.5cm}

Note that the branch flow equations in (1)-(4) are obtained by relaxing the nodal voltage angles as described in \cite{MFarivarLow}. For a radial system, the nodal voltage angle can be exactly obtained from the OPF result. Also, notice that $v_i = |V_i|^2$ and $l_{ij} = |I_{ij}|^2$. It should be noted that (1)-(3) are linear in power flow variables i.e. $P_{ij}, Q_{ij}, v_i \text{ and } l_{ij}$. However, (4) is nonlinear in problem variables (a quadratic equality constraint). 

\subsection{PV Hosting Capacity as D-OPF Problem}
The PV hosting problem is defined as a maximization problem, where the objective is to obtain the maximum PV penetration that can be accommodated by the distribution feeder without violating feeder's operation limits. We formulate the PV hosting capacity problem as a maximization problem with the objective of maximizing the sum of total power injected by the distributed generators connected to prespecified nodes without violating distribution system's operating constraints. We model following constraints in the problem: a) voltage constraints - node voltages should be within prespecified ANSI limits (0.95 pu - 1.05 pu) (8); b) branch thermal limits - line currents should be within the prespecified line thermal ratings (9); c) reverse power flow - no reverse power flow is allowed at the substation bus (10). 


The D-OPF problem to obtain PV hosting capacity of a feeder is defined in (5), subject to power flow constraints in (1)-(4) and additional operational constraints in (6)-(10). Note that the maximization problem can be converted to an equivalent minimization problem by minimizing the negative sum of power injection from all PVs i.e., $-\sum_{i\in \mathcal{N}}{p^{PV}_i}$.  

\begin{equation}
\small
\text{Maximize.} \sum_{i\in \mathcal{N}}{p^{PV}_i}	
\end{equation}
Subject to: (1)-(4), and
\begin{eqnarray}
p_i = p^{PV}_i - p^{load}_i \hspace{1.5 cm} &\forall i \in \mathcal{N}\\
p^{PV}_{i,l} \leq p^{PV}_i \leq p^{PV}_{i,u}\hspace{1.5 cm} &\forall i \in \mathcal{N} \\
(0.95)^2 \leq v_i \leq (1.05)^2  \hspace{1.5 cm} &\forall i \in \mathcal{N}\\
l_{ij} \leq (I_{ij,rated})^2 \hspace{1.5 cm} &\forall \{ij\} \in \mathcal{E}\\
P_{sub} \geq = 0 \hspace{1.5 cm} &
\end{eqnarray}
where, $p^{load}_i$ is the load connected at bus $i$, $p^{PV}_i$ is the rating of PV panel connected to bus $i$, $p^{PV}_{i,l}=0$ and $p^{PV}_{i,u}$ is the maximum PV panel rating that can be connected to bus $i$, $I_{ij,rated}$ is thermal rating of line $\{ij\}$, and $P_{sub}$ is the active power flow out of the substation node. 

\subsection{Convex Relaxation - SOCP Formulation} 
As detailed in several related literature, the quadratic equality constraint (4) in power flow model makes the D-OPF problem for hosting capacity in (5) nonconvex. As such this problem is a quadratically constrained problem and difficult to solve. By relaxing (4) to an inequality as described in (11) a SOCP problem is obtained. 
\begin{equation}
  v_il_{ij} \geq P_{ij}^2 + Q_{ij}^2  \hspace{2.5 cm} \forall i \in \mathcal{N}
\end{equation}

The SOCP-relaxation for D-OPF based PV hosting capacity problem  is detailed as follows. 
\begin{equation}
\small
\text{Maximize.} \sum_{i\in \mathcal{N}}{p^{PV}_i}	
\end{equation}
Subject to: (1)-(3), (6)-(10), and (11).

\subsection{Exactness of the Relaxed D-OPF SOCP Problem}
As described in \cite{Low2}, the SOCP relaxation includes all points inside the second-order cone, however, minimizing a function of current (or power flow) for unconstrained problem yields a feasible solution that lies at the boundary of the cone. Unfortunately, for the case of PV maximization problem, the maximum PV capacity a feeder can accommodate is bounded by upper limits on bus voltages. It has been demonstrated in \cite{Low2} that when upper bounds on node voltages are binding, the SOCP relaxation is not exact. This is because the SOCP solution lies inside the cone and not at the boundary. Due to this SOCP-relaxation of PV hosting capacity problem results in optimal solutions that are not feasible with respect to original quadratic equality constraint. 



\section{Iterative SOCP Formulation for Exact Solution of Relaxed D-OPF Model}
In this section, we propose a convex iteration technique to obtain an optimal and feasible D-OPF solution over multiple SOCP iterations of the relaxed OPF formulation. The proposed approach is inspired from a similar method proposed in \cite{Filter} to design optimal filter parameters using iterative SOCP approach. The solution obtained by solving the iterative procedure will be feasible with respect to the original (before relaxation) D-OPF problem only if 
 the difference between the $P_{ij}^2 + Q_{ij}^2 $ and $v_il_{ij}$ is gradually reduced to zero. 
 
The proposed method is designed to specifically achieve this equality over successive iterations. We define an error term, $e^{(k)}$, measuring the feasibility gap at $k^{th}$ iteration  that is equal to $(P_{ij}^{(k)})^2 + (Q_{ij}^{(k)})^2 - v_i^{(k)} l_{ij}^{(k)}$, where, $P_{ij}^{(k)},Q_{ij}^{(k)},v_{i}^{(k)},l_{ij}^{(k)}$ are power flow variables obtained by solving $k^{th}$ iteration of relaxed D-OPF SOCP model. The objective is to gradually drive the feasibility gap, $e^{(k)}$, to zero over successive iterations. This is achieved by enforcing additional constraint (13), where $\gamma^k < 1$. Notice that $e^{(k)} \leq 0$. Thus, using (13), the feasibility gap, $e^{(k)}$ will gradually increase to zero. Note that,  $\gamma^k$ defines the ratio of error at current iteration and previous iteration and is a tunable parameter (see \cite{Filter} for details). 
\begin{equation}
    \small
    e^{(k)} \geq \gamma^k e^{(k-1)} 
\end{equation}

Also, by substituting for $e^{(k)} = (P_{ij}^{(k)})^2 + (Q_{ij}^{(k)})^2 - v_i^{(k)} l_{ij}^{(k)}$ and using, $\gamma$ where, $\gamma^k \leq \gamma < 1$, the error at $k^{th}$ SOCP-iteration can be expressed using (14).
\begin{equation}
    \small
   e^{(k)} \geq \gamma \times \left((P_{ij}^{(k-1)})^2 + (Q_{ij}^{(k-1)})^2 - v_i^{(k-1)} l_{ij}^{(k-1)}\right)
\end{equation}

As constraint (14) is nonconvex, we linearize it using first-order Taylor series approximation in (15). 
\begin{equation}
\small
\begin{split}
2P_{ij}^{(k-1)}\Delta P_{ij}^{(k)}+2Q_{ij}^{(k-1)}\Delta Q_{ij}^{(k)}-l_{ij}^{(k-1)}\Delta v_i^{(k)}-v_i^{(k-1)}\Delta l_{ij}^{(k)} \geq \\ (\gamma-1)\left((P_{ij}^{(k-1)})^2+(Q_{ij}^{(k-1)})^2-v_i^{(k-1)} l_{ij}^{(k-1)}\right)
\end{split}
\end{equation}
where, at $k^{th}$ iteration, $P_{ij}^{(k-1)}, Q_{ij}^{(k-1)}, v_i^{(k-1)}, \text{ and } l_{ij}^{(k-1)}$ in known from solving $(k-1)^{th}$ SOCP iteration of the problem. Therefore, (15) is linear in unknown, $\Delta P_{ij}^{(k)}, \Delta Q_{ij}^{(k)}, \Delta v_{i}^{(k)} \text{ and } \Delta l_{ij}^{(k)}$. 

Note that in this iterative approach, the actual power flow solution at $k^{th}$ iteration is obtained by updating power flow variables using (16), where $x$ defines set of problem variables. 
\begin{equation}
\small
x^{(k)} = x^{(k-1)} + \alpha \Delta x^{(k)}
\end{equation}
where, the acceleration factor, $0 < \alpha <1 $ and $x^{(k-1)}$ is the variable obtained at previous iteration and change in variables, $\Delta x^{(k)}$ is determined at the current iteration. 

The iterative SOCP model at the $k^{th}$ iteration after including (15) is detailed in (17)-(23).

At $k^{th}$ SOCP-iteration, 
\begin{equation}
\small
\text{Minimize.} -\sum_{i\in \mathcal{N}}{\Delta p^{PV(k)}_i}	
\end{equation}

Subject to: 
\begin{small}
\begin{equation}
\begin{split}
  (P^{(k-1)}_{ij} + \Delta P^{k}_{ij}) = \sum_{k:j\to k}( P^{(k-1)}_{jk} + \Delta P^{k}_{jk} )  + \\ r_{ij}(l^{(k-1)}_{ij} + \Delta l^{(k)}_{ij}) + (p_j^{PV,(k)} + \Delta p_j^{PV,(k)})-p_j^{load} 
\end{split}
\end{equation}

\begin{equation}
\begin{split}
 (Q^{(k-1)}_{ij} + \Delta Q^{k}_{ij}) = \sum_{k:j\to k}( Q^{(k-1)}_{jk} + \Delta Q^{k}_{jk} ) + \\ x_{ij}(l^{(k-1)}_{ij} + \Delta l^{(k)}_{ij}) + q_j^{load}  
\end{split}
\end{equation}

\begin{equation}
\begin{split}
\hspace{-4.0 cm} (v^{(k-1)}_j+\Delta v^{(k)}_j ) = (v^{(k-1)}_i+\Delta v^{(k)}_i) - \\ 2(r_{ij}(P^{(k-1)}_{ij}+\Delta P^{(k)}_{ij})+ x_{ij}(Q^{(k-1)}_{ij}+ \Delta Q^{(k)}_{ij})) + \\
(r_{ij}^2 + x_{ij}^2)(l^{(k-1)}_{ij}+ \Delta l^{(k)}_{ij})
\end{split}
\end{equation}

\begin{equation}
\begin{split}
(v^{(k-1)}_i+ \Delta v^{(k)}_i)(l^{(k-1)}_{ij}+ \Delta l^{(k)}_{ij}) \geq \\((P^{(k-1)}_{ij})^2 + \Delta (P^{(k)}_{ij})^2) + ((Q^{(k-1)}_{ij})^2 + \Delta (Q^{(k)}_{ij})^2)
\end{split}
\end{equation}

\begin{equation}
\begin{split}
\left(2P_{ij}^{(k-1)}\Delta P_{ij}^{(k)} + 2Q_{ij}^{(k-1)}\Delta Q_{ij}^{(k)} -  l_{ij}^{(k-1)}\Delta v_{i}^{(k)} - v_{i}^{(k-1)}\Delta l_{ij}^{(k)}\right) \geq \\ (\gamma-1)\left((P^{(k-1)}_{ij})^2 + (Q^{(k-1)}_{ij})^2 - v^{(k-1)}_{i}l^{(k-1)}_{ij}\right) 
\end{split}
\end{equation}

\begin{equation}
\left(p^{PV}_{i,l} - p_i^{PV,(k)}\right) \leq \Delta p_i^{PV,(k)} \leq \left(p^{PV}_{i,u} -  p_i^{PV,(k)}\right)
\end{equation}
\end{small}
\vspace{-0.2cm}
\begin{equation}
\nonumber \text{and constraints (8)-(10)}
\end{equation}

The decision variables are $\Delta P^{(k)}, \Delta Q^{(k)}, \Delta v^{(k)} , \Delta l^{(k)} $ and $\Delta p_i^{PV,(k)}$. After solving the $kth$ iteration, the variables are updated using (16). The iterative procedure will stop when $|e^{(k)}| \leq \varepsilon$, where $|e^{(k)|}$ is the absolute value of maximum feasibility gap at $k^{th}$ iteration and $\varepsilon$ is a specified tolerance. 

\vspace{-0.2cm}
\section{Results and Discussion}
The proposed iterative method is validated on two test feeders: IEEE 13-bus and IEEE-123 bus \cite{testfeeder}. The test feeders used in the system are first converted into equivalent positive sequence model using OpenDSS. All the simulations are done using MATLAB. The initial conditions are obtained from the solution of an equivalent linearized D-OPF problem solved using CPLEX 12.7 \cite{Low1}. The relaxed-SOCP and iterative-SOCP problems are formulated in MATLAB environment and solved using fmincon function in MATLAB optimization toolbox. The exact nonlinear OPF problem is also formulated and solved for validation of the results. The PV hosting capacity problem is solved for maximum and minimum load conditions. For maximum load, all loads are at their peak load, while for minimum load a load multiplier of 0.3 is assumed. 

\vspace{-0.2cm}
\subsection{Feasibility of  relaxed problem}
In this section, we demonstrate that for PV maximization problem, the relaxed SOCP model leads to feasibility violations. By relaxing the constraints, there is an increase in the feasible solution space of the problem. The solution obtained after solving the relaxed-SOCP model are feasible when the \textit{quadratic  equality} constraints are satisfied. Therefore, the SOCP relaxation is exact if the obtained solution lies on the surface of the cone i.e. the solution to relaxed problem indeed satisfies $v_il_{ij} = P_{ij}^2 + Q_{ij}^2$.   

For both IEEE-13 node and IEEE-123 node test systems, we solve the relaxed-SOCP model to obtain the PV hosting capacity. The feasibility violation is measured by computing $e_{ij} = P_{ij}^2 + Q_{ij}^2 - v_il_{ij}$ for all $\{ij\} \in \mathcal{E}$. The plots for feasibility violation represented using $e_{ij}$ for each branch corresponding to IEEE-13 bus and IEEE-123 bus systems are shown in Fig.3. Note that with relaxed constraints, the maximum feasibility violation for IEEE-13 bus system is -87.71 and -71.44 at the minimum and maximum loading conditions respectively. Similarly, for the IEEE-123 node system, the maximum feasibility violation is -39.14 and -38.96 at the minimum and maximum loading condition, respectively. It can be concluded that the relaxed quadratic constraints lead to a solution that is not exact wrt. the original quadratic equality constraints or the original nonlinear-OPF model. 


\begin{figure}[t]
\centering
\includegraphics[trim={1cm 0 2cm 0cm},clip, width=0.95\linewidth]{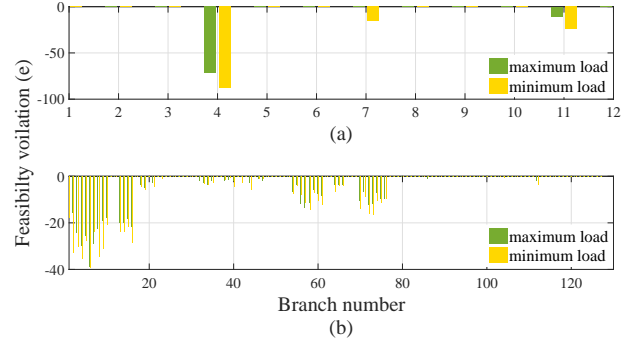}
\vspace{-0.2cm}
	\caption{feasibility violation in a) IEEE 13 bus  b) IEEE 123 bus system }
	\vspace{-0.4cm}
	\label{fig:3}
	\vspace{-0.0cm}
	\end{figure}

\subsection{Iterative SOCP Model for PV Hosting Capacity Problem}
It is clear from the previous section that the solutions obtained from SOCP relaxation of D-OPF model for PV hosting problem are not exact wrt. the original nonlinear-OPF model. In this section, we employ the proposed iterative framework to obtain the exact OPF solutions by solving multiple iterations of the relaxed-SOCP problem. The results for PV hosting capacity are shown for both IEEE-13 and IEEE-123 systems.

\subsubsection{IEEE-13 bus system}
The IEEE-13 bus is a small but heavily loaded feeder with the total demand of 1155 kW operating at 4.16 kV. The objective is to maximize the power supplied by the distributed PVs in the distribution system while satisfying the operating constraints detailed in Section II.B.

For the optimization, it is assumed that PVs can be placed at all the nodes in the distribution feeder. For this system the maximum rating of each PV system is assumed to be 400 kW. The PV hosting capacity at minimum and maximum loading conditions are shown in Table I. Note that the hosting capacity obtained from the proposed iterative method is equal to the hosting capacity obtained by solving the actual nonlinear-OPF model (one with nonlinear quadratic equality constraints). Therefore, this case validates the proposed convex iteration method leads to both feasible and optimal solution for the original nonconvex OPF problem. In this case, the hosting capacity is limited by the reverse power flow condition. 

Now we turn our attention to the satisfaction of equality constraints over successive relaxed-SOCP iterations. To analyze the feasibility gap at each iteration and to show that indeed it decreases to the desired error bound, we measure $e_{ij} = P_{ij}^2 + Q_{ij}^2 - v_il_{ij}$ for all $\{ij\} \in \mathcal{E}$ and identify its maximum absolute value at each iteration. The plot for maximum feasibility error observed at each iteration is shown in Fig. 4. Note that the maximum feasibility error at first iteration is 0.0071 and it decreases to 0.00061 at the $12th$ iteration. Similarly, for the maximum loading condition, the error at first iteration is 0.0367 and it decreased to 0.0032 at  $12th$ iteration. It should be noted that we do not need to simulate 12 iterations. Since the feasibility error is small to begin with, due to introduction of the additional linear constraints, the iterations can be stopped sooner.  

	\begin{figure}[t]
\centering
\includegraphics[width=3.7in]{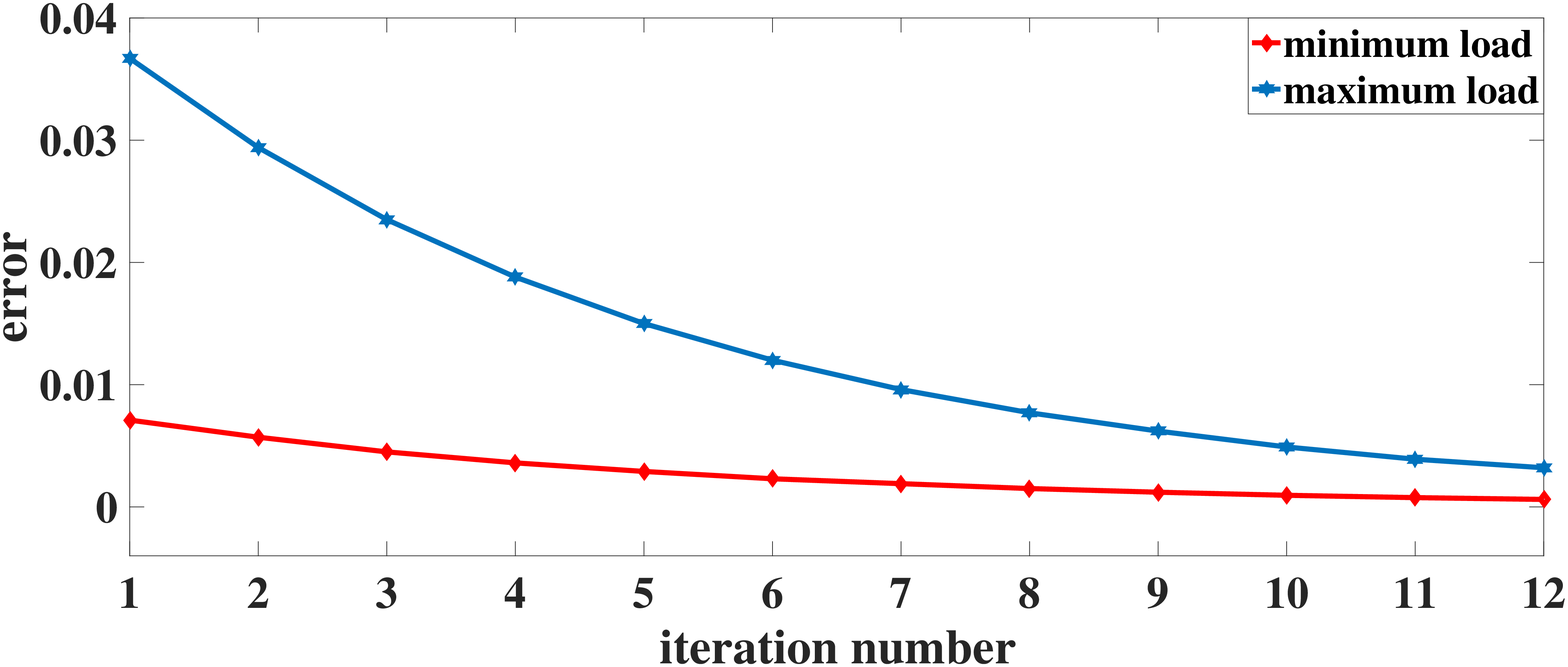}
\vspace{-0.8cm}
\caption{IEEE-13 bus:Reduction in feasibility gap vs. number of iterations}
	\vspace{-0.5cm}
	\label{fig:7}
	\end{figure}
	
\begin{figure}[t]
\centering
\includegraphics[width=3.7in]{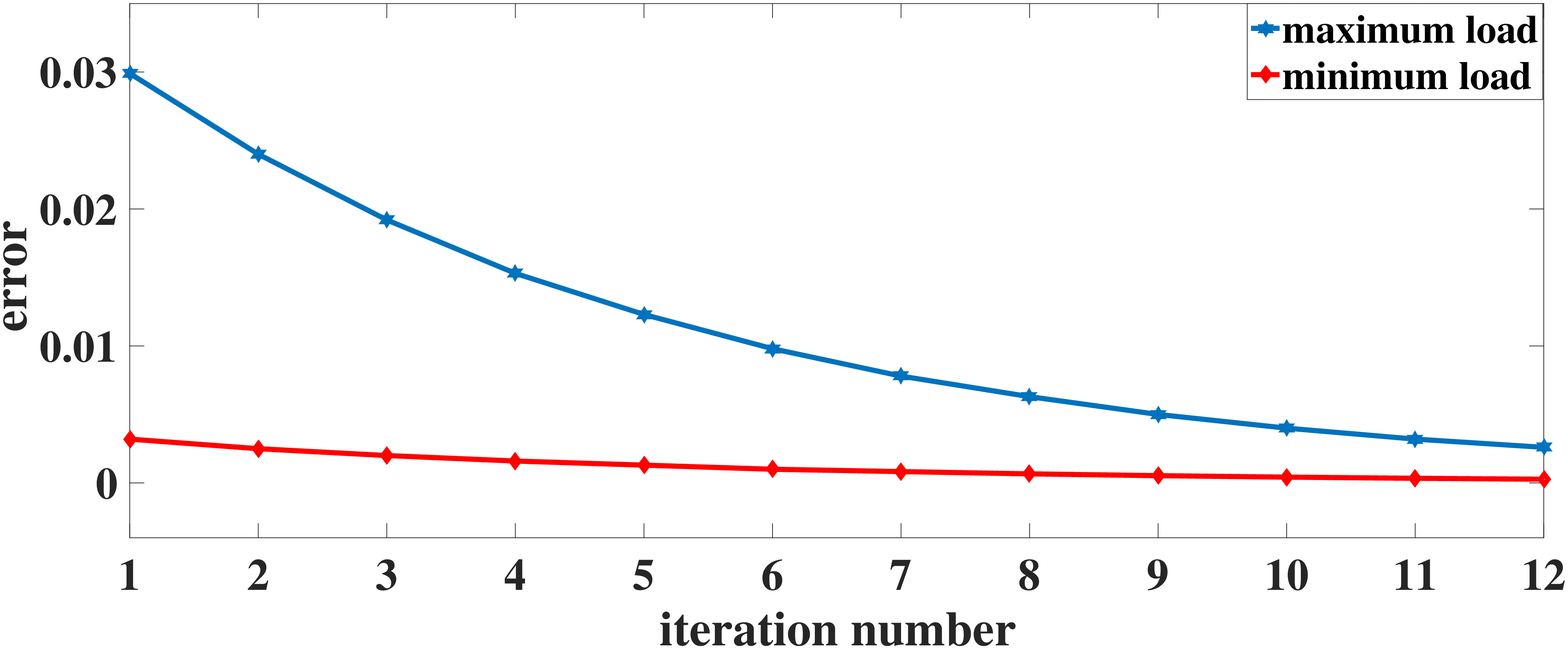}
\vspace{-0.5cm}
\caption{IEEE-123 bus:Reduction in feasibility gap vs. number of iterations}
	\vspace{-1.2cm}
	\label{fig:10}
	\vspace{0.9cm}
\end{figure}	
	

\begin{table}[t]
\centering
\caption{PV Hosting capacity for IEEE-13 and 123 bus system (in kW)}
\vspace{-0.1cm}
\begin{tabular}{c|c|c|c}
\hline
\hline
Test Feeder & loading condition & Nonlinear OPF & Iterative OPF\\
\hline
\multirow{2}{*}{IEEE-13 bus}& minimum & 350.6& 350.5\\
\cline{2-4}
& maximum&1185.7&1185.9\\
\hline
\multirow{2}{*}{IEEE-123 bus}& minimum & 350.42& 350.45\\
\cline{2-4}
& maximum&1184.1&1184.6\\
\hline
\hline
\end{tabular}

	\vspace{-0.4cm}
\end{table}


\subsection{IEEE -123 bus system}	
The IEEE-123 bus is a larger feeder with the total load demand of 1163 kW operating at 2.4 kV. 
It is assumed that PVs can be placed at all the nodes in the distribution system and the maximum rating of each PV is bounded by 50 kW. The results obtained for the PV hosting capacity at the minimum and maximum loading conditions are shown in Table I. 

Next, we analyze the feasibility gap at each iteration of the proposed iterative SOCP formulation. Due to additional constraints, the maximum error for the system at the minimum load condition at first iteration is 0.0032 decreases over successive iterations. Similarly, for the maximum load condition, the error at first iteration is 0.0311 to 0.0032 at $12^{th}$ iteration. It is to be noted that as the number of iteration increases the feasibility gap is driven towards zero. We also calculate the PV hosting capacity using original nonlinear-OPF model for both loading conditions (see Table 1). The results show that the proposed iterative approach leads to same solution and hence is successful in achieving an optimal and feasible solution for the PV hosting problem. 





 





\section{Conclusions}

The optimal power flow problems are nonconvex and difficult to solve. Although several relaxation methods have been proposed, the exactness of the solutions obtained from the relaxed models remain a concern. In this paper, we identify one such problem related to second-order cone relaxation of branch-power flow based OPF model when attempting to solve a maximization problem. The problem of maximizing the feeder's PV hosting capacity is modeled as an OPF problem and relaxed using SOCP relaxation. On solving the relaxed model, it is demonstrated that the results are not feasible with respect to the original nonlinear problem. In order to solve this problem, we have proposed a convex iteration technique that by solving multiple iterations of the relaxed-SOCP problem results in an exact solution for the original nonlinear-OPF problem. The proposed approach is validated on IEEE-13 and IEEE-123 bus test feeders. The results demonstrate that the proposed iterative method reduces the feasibility gap observed due to SOCP relaxations over multiple convex iterations of the relaxed-SOCP problem.




\bibliographystyle{IEEEtran}
%
%
  
\bibliography{references}

\begin{thebibliography}{10}
\providecommand{\url}[1]{#1}
\csname url@samestyle\endcsname
\providecommand{\newblock}{\relax}
\providecommand{\bibinfo}[2]{#2}
\providecommand{\BIBentrySTDinterwordspacing}{\spaceskip=0pt\relax}
\providecommand{\BIBentryALTinterwordstretchfactor}{4}
\providecommand{\BIBentryALTinterwordspacing}{\spaceskip=\fontdimen2\font plus
\BIBentryALTinterwordstretchfactor\fontdimen3\font minus
  \fontdimen4\font\relax}
\providecommand{\BIBforeignlanguage}[2]{{%
\expandafter\ifx\csname l@#1\endcsname\relax
\typeout{** WARNING: IEEEtran.bst: No hyphenation pattern has been}%
\typeout{** loaded for the language `#1'. Using the pattern for}%
\typeout{** the default language instead.}%
\else
\language=\csname l@#1\endcsname
\fi
#2}}
\providecommand{\BIBdecl}{\relax}
\BIBdecl

\bibitem{OPFsurvey}
W.~Sheng, K.~Liu, and S.~Cheng, ``Optimal power flow algorithm and analysis in
  distribution system considering distributed generation,'' \emph{IET
  Generation, Transmission Distribution}, vol.~8, no.~2, pp. 261--272, February
  2014.

\bibitem{OPFsurvey2}
M.~Huneault and F.~D. Galiana, ``A survey of the optimal power flow
  literature,'' \emph{IEEE Transactions on Power Systems}, vol.~6, no.~2, pp.
  762--770, May 1991.

\bibitem{SHLowCR1}
S.~H. Low, ``Convex relaxation of optimal power flow—part i: Formulations and
  equivalence,'' \emph{IEEE Transactions on Control of Network Systems},
  vol.~1, no.~1, pp. 15--27, March 2014.

\bibitem{approx}
Z.~Yuan and M.~R. Hesamzadeh, ``Improving the accuracy of second-order cone ac
  optimal power flow by convex approximations,'' in \emph{2018 IEEE Innovative
  Smart Grid Technologies - Asia (ISGT Asia)}, May 2018, pp. 172--177.

\bibitem{XBai}
X.~Bai and H.~Wei, ``Semi-definite programming-based method for
  security-constrained unit commitment with operational and optimal power flow
  constraints,'' \emph{IET Generation, Transmission Distribution}, vol.~3,
  no.~2, pp. 182--197, February 2009.

\bibitem{RAJabr}
R.~A. Jabr, ``Radial distribution load flow using conic programming,''
  \emph{IEEE Transactions on Power Systems}, vol.~21, no.~3, pp. 1458--1459,
  Aug 2006.

\bibitem{DKMolzahn}
D.~K. Molzahn and I.~A. Hiskens, ``Convex relaxations of optimal power flow
  problems: An illustrative example,'' \emph{IEEE Transactions on Circuits and
  Systems I: Regular Papers}, vol.~63, no.~5, pp. 650--660, May 2016.

\bibitem{Low1}
M.~Farivar and S.~H. Low, ``Branch flow model: Relaxations and convexification:
  {Part I},'' \emph{IEEE Transactions on Power Systems}, vol.~28, no.~3, pp.
  2554--2564, Aug 2013.

\bibitem{Low2}
S.~H. Low, ``Convex relaxation of optimal power flow {Part II}: Exactness,''
  \emph{IEEE Transactions on Control of Network Systems}, vol.~1, no.~2, pp.
  177--189, June 2014.

\bibitem{MFarivar}
M.~Farivar, R.~Neal, C.~Clarke, and S.~Low, ``Optimal inverter var control in
  distribution systems with high pv penetration,'' in \emph{2012 IEEE Power and
  Energy Society General Meeting}, July 2012, pp. 1--7.

\bibitem{LGan}
L.~Gan, N.~Li, U.~Topcu, and S.~H. Low, ``Exact convex relaxation of optimal
  power flow in radial networks,'' \emph{IEEE Transactions on Automatic
  Control}, vol.~60, no.~1, pp. 72--87, Jan 2015.

\bibitem{Sbose}
S.~Bose, S.~H. Low, T.~Teeraratkul, and B.~Hassibi, ``Equivalent relaxations of
  optimal power flow,'' \emph{IEEE Transactions on Automatic Control}, vol.~60,
  no.~3, pp. 729--742, March 2015.

\bibitem{iter2}
E.~Grover-Silva, R.~Girard, and G.~Kariniotakis, ``Multi-temporal optimal power
  flow for assessing the renewable generation hosting capacity of an active
  distribution system,'' in \emph{2016 IEEE/PES Transmission and Distribution
  Conference and Exposition (T D)}, May 2016, pp. 1--5.

\bibitem{iter1}
S.~Y. Abdelouadoud, R.~Girard, F.~Neirac, and T.~Guiot, ``Iterative linear cuts
  strenghtening the second-order cone relaxation of the distribution system
  optimal power flow problem,'' in \emph{2014 IEEE PES T D Conference and
  Exposition}, April 2014, pp. 1--4.

\bibitem{Filter}
A.~Jiang and H.~K. Kwan, ``Minimax design of iir digital filters using
  iterative socp,'' \emph{IEEE Transactions on Circuits and Systems I: Regular
  Papers}, vol.~57, no.~6, pp. 1326--1337, June 2010.

\bibitem{MFarivarLow}
M.~Farivar and S.~H. Low, ``Branch flow model: Relaxations and
  convexification—part i,'' \emph{IEEE Transactions on Power Systems},
  vol.~28, no.~3, pp. 2554--2564, Aug 2013.

\bibitem{testfeeder}
K.~P. Schneider, B.~A. Mather, B.~C. Pal, C.~. Ten, G.~J. Shirek, H.~Zhu, J.~C.
  Fuller, J.~L.~R. Pereira, L.~F. Ochoa, L.~R. de~Araujo, R.~C. Dugan,
  S.~Matthias, S.~Paudyal, T.~E. McDermott, and W.~Kersting, ``Analytic
  considerations and design basis for the ieee distribution test feeders,''
  \emph{IEEE Transactions on Power Systems}, vol.~33, no.~3, pp. 3181--3188,
  May 2018.

\end{thebibliography}


\end{document}